\documentstyle{article}
\begin{document}
\baselineskip+6pt \small \begin{center}{\textit{ In the name of
Allah, the Beneficent, the Merciful}}\end{center}
\large

\begin{center}{ON THE FIELD OF DIFFERENTIAL RATIONAL INVARIANTS OF A SUBGROUP
OF AFFINE GROUP(ORDINARY DIFFERENTIAL CASE)} \end{center}

\begin{center} {Ural Bekbaev \footnote[1]{e-mail: bekbaev@science.upm.edu.my}}
\end{center}

\begin{center} {Department of Mathematics $\&$ Institute for Mathematical Research,}\end{center}
\begin{center} {FS, UPM, 43400, Serdang, Selangor, Malaysia.}\end{center}

\begin{abstract}{An ordinary differential field  $(F,d)$ of characteristic zero, a
subgroup $H$ of affine group $ GL(n,C)\propto C^n$ with respect to
its identical representation in $F^n$ and the following two fields
of differential rational functions in $x=(x_1,x_2,...,x_n)$-column
vector,
$$C\langle x, d\rangle^H=\{f^d\langle x\rangle \in C\langle x, d\rangle :
f^d\langle hx+ h_0\rangle = f^d\langle x\rangle\ \mbox{for any} \
(h,h_0)\in H \},$$
$$C\langle x, d\rangle^{(F^*,H)}=\{f^d\langle x\rangle \in C\langle x, d\rangle :
f^{g^{-1}d}\langle hx+ h_0\rangle = f^d\langle x\rangle\ \mbox{for
any} \ g\in F^* \ \mbox{and} \ (h,h_0)\in H \}$$ are considered,
where $C$ is the constant field of $(F,d)$ and $C\langle x,
d\rangle$ is the field of differential rational functions in
$x_1,x_2,...,x_n$ over $C$. The field $C\langle x, d\rangle^H$
($C\langle x, d\rangle^{(F^*,H)}$) is an important tool in the
equivalence problem of paths( respect. curves) in Differential
Geometry with respect to the motion group $H$. In this paper an
pure algebraic approach is offered to describe these fields. The
field $C\langle x, d\rangle^{(F^*,H)}$ and its relation with
$C\langle x, d\rangle^H$ are investigated. It is shown also that
$C\langle x, d\rangle^H$ can be derived from some algebraic (
without derivatives) invariants of $H$. \vspace{0.1cm}
} \end{abstract}
\vspace{0.7cm}

{\bf  1. Introduction}

Let $n>1$ be a natural number, $H$ be a subgroup of the affine
group $ GL(n,R)\propto R^n $, and $G= Diff(B)$ be the group of
diffeomorphisms of the open interval $B=(0,1)$. Two infinitely
smooth paths $u:B\rightarrow R^n$, $v:B\rightarrow R^n$ are said
to be $(G,H)$-equivalent if there are such $t \in G$ and
$(h,h_0)\in H$ that $ v(s)=hu(t(s))+h_0 \quad \mbox{for any}\quad
s\in B.$

A function $f^d(u(t))$ of $u(t)=(u_1(t),. . ., u_n(t))$ and its
finite number of derivatives relative to $d= \frac{d}{dt}\quad$ is
said to be invariant(more exactly, $ (G,H)$- invariant) if the
equality
$$f^d(u(t(s)))= f^{\delta}(hu(t(s)) + h_0)$$ is valid for any $t\in G$  and $(h,h_o)\in H$,where $\delta= \frac{d}{ds}$.

In the terms of inverse function $s=s(t)$ the expressions
$\frac{du(t(s))}{dt(s)}$ and $\frac{du(t(s))}{ds}$ can be
rewritten as $\frac{du(t)}{dt}$ and
$\frac{du(t)}{ds(t)}=\frac{1}{s'(t)}\frac{du(t)}{dt}$,
respectively. Therefore the above invariantness of $f^d(u(t))$ can
be written in the form $f^d(u(t))= f^{\delta}(hu(t) + h_0)$, where
$\delta= s'(t)^{-1}\frac{d}{dt}$.

Let $t$ run $B$, $F= C^{\infty}(B)$, and $(F,d)$ be the
differential ring of infinitely smooth functions on $B$ relative
to differential operator $\frac{d}{dt}$. The constant subring,i.e.
$\{ a\in F: d(a)=0\}$, of $(F,d)$ is $R$. Now every infinitely
smooth path $u:B\rightarrow R^n$ can be considered as an element
of differential module $(F^n;d )$, where $d$ acts on elements of
$F^n$ coordinate-wisely. If elements of this module are considered
as column vectors the above transformations of $u$ and $d$ look
like $$u=(u_1,. . .,u_n)\mapsto hu+ h_0,\quad d \mapsto g^{-1}d$$,
where  $g$ is some invertible element of $F$ and $(h,h_0)\in H $.

Therefore the following algebraic approach to the above invariants
is natural. Let $(F,d)$ be any differential field i.e. $F$ is a
field with a fixed $d:F \rightarrow F$ for which
$$ a)\quad d(a+b)= d(a)+d(b),\quad\quad  b)\quad d(ab)=d(a)b+ad(b)
\quad\mbox{for any}\quad a,b\in F.$$
Let $C$ stand for the constant subfield of $(F,d)$ i.e.
  $C= \{ a\in F: d(a)=0\}$
, $H\subset GL(n,C)\propto C^n$ be a subgroup, $F^*=F\setminus
\{0\}$.

One can consider the following $(F^*,H)$-equivalence of pairs
$(u,d_u)$, $(v, d_v)$, where $u,v \in F^n$ and $d_u:F \rightarrow
F $, $d_v:F \rightarrow F $ are differential operators, as an
algebraic analogue of the above mentioned equivalence of paths.

{\bf Definition 1.} \textit{ Pairs $(u,d_u)$, $(v, d_v)$ are said
to be $(F^*,H)$-equivalent if $v =hu+ h_0, d_v= g^{-1}d_u $ for
some $g\in F^*$ and $(h,h_0)\in H$.}

  Let in future $x_1,. . .,x_n$ be $d$-differential algebraic independent
 variables over $F$ and
$x$ stand for the column vector with coordinates $x_1,. . .,x_n$.
We use the following notations: $C[x]$- the ring of polynomials in
variables $x_1,. . .,x_n$ over $C$; $C(x)$- the field of rational
functions in variables $x_1,. . .,x_n$ over $C$;$C\{x,d\}$- the
ring of differential polynomials in $x_1,. . .,x_n$ over $C$ i.e.
$$C\{x,d \}=C[x,dx,d^2x,...,d^mx,...]$$, where
$d^mx=(d^mx_1,...,d^mx_n)$ and  $C\langle x, d\rangle$ the field
of differential rational functions in $x$ over $C$.

{\bf Definition 2.} \textit{ An element $f^d\langle x\rangle \in
C\langle x, d\rangle$ is said to be  $(F^*,H)$ invariant \\($F^*$
invariant; $H$ invariant)-
 if the equality $$f^{g{-1}d}\langle hx+ h_0\rangle =
f^d\langle x\rangle $$ (respect.$f^{g{-1}d}\langle x\rangle =
f^d\langle x\rangle$;$f^d\langle hx+ h_0\rangle = f^d\langle
x\rangle$) is valid for any $g\in F^*$ and $(h,h_0)\in H$.}

Let us denote the set of all $(F^*,H)$ (respect. $F^*$; $H$)-
invariant differential rational functions over $C$ by $C\langle x,
d\rangle^{(F^*,H)}$ (respect. $C\langle x, d\rangle^{F^*}$;
$C\langle x, d\rangle^H$).

Importance of differential invariants ( for example such as,
curvature, torsion) is out of the question. Therefore
investigation of the field $C\langle x, d\rangle^{(F^*,H)}$, as
far as it is an algebraic analog of such differential invariants,
is of interest. In Differential Geometry usually geometric
methods, for example Cartan's moving frame method, are used to
obtain differential invariants of curves with respect to a motion
group $H$.  It is clear that the transformation $d \mapsto
g^{-1}d$ corresponds to the ordinary change of parameter in theory
of curves. In this paper we are going to offer a pure algebraic
approach to describe the field $C\langle x, d\rangle^{(F^*,H)}$.
Due to the algebraic character of our approach the considered case
may go even beyond the current needs of Differential Geometry. The
finite group case of $H$ is considered in [1].

As to the field $C\langle x, d\rangle^H$ it is an algebraic analog
of differential invariants of paths with respect to the motion
group $H$. This field is investigated for all classical subgroups
 of the affine group $GL(n,R)\propto R^n$. The corresponding results
 and references can be found in [2].

 The results of this paper are presented in [3,4].

{\bf 2. Preliminary }

In future $(F;d)$ stands for an ordinary differential field of
characteristic zero and $C$ is its constant field. It is assumed
that the differential field $(F,d)$ holds the following property:

If $q^d\{y\}\in F\{y, d \}$ is such a polynomial that $q^d\{a\}=0
$ for any $a\in F$ then $q^d\{y\}$ is a zero polynomial.

This condition is equivalent to $F\neq C$ due to [5, p.139].

{\bf Proposition 1.}\textit{Every non constant $f^{d}\langle x
\rangle \in F\langle x, d \rangle$ is $d$-algebraic independent
over $F$.}

{\bf Proof.} Represent non constant $f^{d}\langle x \rangle $ as
an irreducible ratio $ \frac{P\{x\}}{Q\{x\}}$, where
$P\{x\},Q\{x\} \in F\{x, d \}$. Assume that, for example, $x_1$
occurs in $P\{x\}$ or $Q\{x\}$ and define $d$-order of $f^d\langle
x \rangle $ with respect to $x_1$ as the maximal $k$ that $d^kx_1$
occurs in $P\{x\}$ or $Q\{x\}$.

To show $d$-algebraic independence of $f^d\langle x \rangle $ over
$F$ it is enough to show that the similar order of $df^{d}\langle
x \rangle $ is $k+1$. Indeed the coefficient at $d^{k+1}x_1$ in
the numerator of $d\frac{P\{x\}}{Q\{x\}}$ is equal to
$$\frac{\partial P\{x\}}{\partial d^kx_1}Q\{x\}-\frac{\partial
 Q\{x\}}{\partial d^kx_1}P\{x\}$$
which is nonzero because characteristic of $F$ is zero and the
ring $F\{x,d\}$ is a factorial ring. Therefore for the $d$-order
of $df^{d}\langle x \rangle $ with respect to $x_1$ one has $k+1$.

Now  let us show that if $H$ is a subgroup then one can find such
a nonzero $p^d\langle x\rangle \in C\langle x,d\rangle$ for which
\begin{eqnarray}\begin{array}{c}p^{g^{-1}d}\langle hx+h_0\rangle =
g^{-1}p^d\langle x\rangle, \end{array}\end{eqnarray}
for any $g\in F^*,\quad (h,h_0)\in H$.

Here are two examples of nonzero $p^d\langle x\rangle$ for which
property (1) holds in the case of $H= GL(n,C)\propto C^n$ and
therefore equality (1) is valid when $H$ is any subgroup of the
affine group as well.

{ \bf Example 1. } Let $g\in F^*$ and $\delta$ stand for
$g^{-1}d$. The following relation can be verified.
$$d^k= \sum_{i=1}^k \Phi^d_{k,i}\{g\}\delta^i$$, where $\Phi^d_{k,i}\{g\}=
\sum C_{\alpha}g^{[\alpha ]}$, $\alpha=
(\alpha_1,\alpha_2,...,\alpha_k)$, $ g^{[\alpha
]}=g^{\alpha_1}(dg)^{\alpha_2}...(d^{k-1}g)^{\alpha_k}$,
$C_{\alpha}=\frac{k!}{\alpha_1!\alpha_2!...\alpha_{k}!\prod_{j=1}^k
(j!)^{\alpha_j}}$, the sum $ \sum $ is taken over all $\alpha$
with nonnegative integer components for which   $\vert \alpha
\vert= \sum_{j=1}^k\alpha_j = i$ and $
1\alpha_1+2\alpha_2+...+k\alpha_k=k$.

In future let us use the notations:$W^d= W^d\{dx\}=
W^d\{dx_1,dx_2,...,dx_n\}= \det[dx,d^2x,..., d^{n}x]$. Let
$W^d_i=W^d_i\{dx,d^{n+1}x\}$, $i=\overline{1,n}$ stand for the
determinant of the matrix obtained from \\ $[dx,d^2x,..., d^nx,
d^{n+1}x]$ by deleting $d^ix$.

 Consider one more variable $y$. It is easy to check that
 $$W^{\delta}(\delta x_1,\delta x_2,...,\delta x_n,\delta y )=g^{-\frac{(n+1)(n+2)}{2}}W^d(dx_1,dx_2,...,dx_n,dy ).$$ Therefore
$\sum^{n+1}_{i=1}(-1)^{n+1-i}W^{\delta}_i\delta^iy=
g^{-\frac{(n+1)(n+2)}{2}}\sum^{n+1}_{i=1}(-1)^{n+1-i}W^d_id^iy=$
$$g^{-\frac{(n+1)(n+2)}{2}}\sum^{n+1}_{i=1}(-1)^{n+1-i}W^d_i\sum^i_{j=1} \Phi^d_{i,j}\{g\}\delta^jy= g^{-\frac{(n+1)(n+2)}{2}}\sum^{n+1}_{j=1}(\sum^{n+1}_{i=j}(-1)^{n+1-i} \Phi^d_{i,j}\{g\}W^d_i)\delta^jy.$$

It implies that for any $j=\overline{1,n+1}$ one has
$W^{\delta}_j=g^{-\frac{(n+1)(n+2)}{2}}\sum^{n+1}_{i=j}(-1)^{i-j}\Phi^d_{i,j}\{g\}W^d_i.$

In particular, if $j= n+1$ then
$W^{\delta}=g^{-\frac{n(n+1)}{2}}W^d,$\\
if $j= n$ then
$W^{\delta}_n=g^{-\frac{n(n+1)}{2}-1}(W^d_{n}-\frac{n(n+1)}{2}\frac{dg}{g}W^d),$\\
if $j= n-1$ then $W^{\delta}_{n-1}= $
$$g^{-\frac{n(n+1)}{2}-2}(W^d_{n-1}-
\frac{n(n-1)}{2}\frac{dg}{g}W^d_{n}+(\frac{(n-1)n(n+1)}{6}d(\frac{dg}{g})+
\frac{(n-1)n(n+1)(3n-2)}{24}(\frac{dg}{g})^2)W^d),$$ if $j= n-2$
then
$$W^{\delta}_{n-2}=g^{-\frac{n(n+1)}{2}-3}(W^d_{n-2}-
\frac{(n-1)(n-2)}{2}\frac{dg}{g}W^d_{n-1}+(\frac{(n-2)(n-1)n}{6}d(\frac{dg}{g})+
$$
$$\frac{(n-2)(n-1)n(3n-5)}{24}(\frac{dg}{g})^2)W^d_n-(\frac{(n+1)!}{24(n-3)!}d^2(\frac{dg}{g})+$$
$$(2n-3)\frac{(n+1)!}{24(n-3)!}\frac{dg}{g}d(\frac{dg}{g})+(n-1)(n-2)\frac{(n+1)!}{48(n-3)!}(\frac{dg}{g})^3)W^d).$$

Therefore
\begin{eqnarray}\begin{array}{c}
\frac{W^{\delta}_n}{W^{\delta}}=g^{-1}(\frac{W^d_n}{W^d}-\frac{n(n+1)}{2}\frac{dg}{g})\end{array}\end{eqnarray}
\begin{eqnarray}\begin{array}{c}
\frac{W^{\delta}_{n-1}}{W^{\delta}}=g^{-2}(\frac{W^d_{n-1}}{W^d}-\frac{n(n-1)}{2}\frac{W^d_{n}}{W^d}\frac{dg}{g}+
 (\frac{(n-1)n(n+1)}{6}d(\frac{dg}{g})+ \frac{(n-1)n(n+1)(3n-2)}{24}(\frac{dg}{g})^2))\end{array}\end{eqnarray}
Due to (2) one has the following two equalities
\begin{eqnarray}\begin{array}{c}
(\frac{W^{\delta}_n}{W^{\delta}})^2=g^{-2}((\frac{W^d_n}{W^d})^2-n(n+1)(\frac{W^d_n}{W^d}\frac{dg}{g}-\frac{n(n+1)}{4}(\frac{dg}{g})^2))\end{array}\end{eqnarray}
$$\delta(\frac{W^{\delta}_n}{W^{\delta}})=g^{-2}(d(\frac{W^d_n}{W^d})-\frac{W^d_n}{W^d}\frac{dg}{g} -\frac{n(n+1)}{2}\frac{d^2g}{g}+ n(n+1)(\frac{dg}{g})^2)$$
 The last equality and (3) imply that
$$\frac{W^{\delta}_{n-1}}{W^{\delta}}+\frac{n-1}{3}\delta(\frac{W^{\delta}_n}{W^{\delta}})=
g^{-2}(\frac{W^d_{n-1}}{W^d}+\frac{n-1}{3}d(\frac{W^d_n}{W^d})\\
-\frac{(n-1)(3n+2)}{6}(\frac{W^d_{n}}{W^d}\frac{dg}{g}-
\frac{n(n+1)}{4}(\frac{dg}{g})^2)).$$ Now one can use this
equality with (4) to get
$$\frac{W^{\delta}_{n-1}}{W^{\delta}}+\frac{n-1}{3}\delta(\frac{W^{\delta}_n}{W^{\delta}})-\frac{(n-1)(3n+2)}{6n(n+1)}(\frac{W^{\delta}_n}{W^{\delta}})^2=
g^{-2}(\frac{W^d_{n-1}}{W^d}+\frac{n-1}{3}d(\frac{W^d_n}{W^d})-\frac{(n-1)(3n+2)}{6n(n+1)}(\frac{W^d_n}{W^d})^2).$$
So for the $p^d_1\langle x\rangle =
\frac{W^d_{n-1}}{W^d}+\frac{n-1}{3}d(\frac{W^d_n}{W^d})-\frac{(n-1)(3n+2)}{6n(n+1)}(\frac{W^d_n}{W^d})^2$
one has $p^{\delta}_1\langle x\rangle =g^{-2}p^d_1\langle
x\rangle$. Therefore $\frac{\delta p^{\delta}_1\langle
x\rangle}{p^{\delta}_1\langle x\rangle}=
g^{-1}(\frac{dp^d_1\langle x\rangle}{p^d_1\langle
x\rangle}-2\frac{dg}{g}).$ Combining it with (2)one can see that
the following nonzero function $p^d\langle x\rangle =
\frac{n(n+1)}{2}\frac{dp^d_1\langle x\rangle}{p^d_1\langle
x\rangle}- 2 \frac{W^d_n}{W^d}$ meets the needed requirements.

{ \bf Example 2.} Note that if one substitutes
$\frac{2}{n(n+1)}\frac{W^d_n}{W^d}$ for $\frac{dg}{g}$ in the
expression
$$\frac{W^d_{n-1}}{W^d}-\frac{n(n-1)}{2}\frac{W^d_{n}}{W^d}\frac{dg}{g}+
 (\frac{(n-1)n(n+1)}{6}d(\frac{dg}{g})+ \frac{(n-1)n(n+1)(3n-2)}{24}(\frac{dg}{g})^2),$$
which is a part of equality (3), he gets the above function
$p^d_1\langle x\rangle $.

In similar way if one substitutes
$\frac{2}{n(n+1)}\frac{W^d_n}{W^d}$ for $\frac{dg}{g}$ in the
corresponding part of the expression for
$\frac{W^{\delta}_{n-2}}{W^{\delta}}$, namely
$$\frac{W^d_{n-2}}{W^d}-\frac{(n-1)(n-2)}{2}\frac{W^d_{n-1}}{W^d}\frac{dg}{g}+
 (\frac{(n-2)(n-1)n}{6}d(\frac{dg}{g})+ \frac{(n-2)(n-1)n(3n-5)}{24}(\frac{dg}{g})^2)\frac{W^d_{n}}{W^d}-$$
$$(\frac{(n+1)!}{24(n-3)!}d^2(\frac{dg}{g})+ (2n-3)\frac{(n+1)!}{24(n-3)!}\frac{dg}{g}d(\frac{dg}{g})+(n-1)(n-2)\frac{(n+1)!}{48(n-3)!}(\frac{dg}{g})^3)),$$
 then he gets a function $p^d_2\langle x\rangle$, for which
$p^{g^{-1}d}_2\langle hx+h_0\rangle = g^{-3}p^d_2\langle x\rangle$
for any $g\in F^*$ and $(h,h_0)\in H.$ Therefore the function
$p^d\langle x\rangle = \frac{p^d_2\langle x\rangle}{p^d_1\langle
x\rangle}$ also holds the needed property (1).

To have some other examples of such $p^d\langle x\rangle$ one can
see [6].

{\bf 3. On $(F^*,H)$- invariants }

Due to the above examples for a given subgroup $H$ of the affine
group $GL(n,C)\propto C^n$ one can consider $\delta (x,d)= \delta=
p^d\langle x\rangle^{-1}d$, where $p^d\langle x\rangle \in
C\langle x,d\rangle$ is a fixed nonzero element holding property
(1). It is evident that the field $C\langle x,d\rangle^{(F^*,H)}$
is invariant with respect to the differential operator $\delta$ ,
so one can consider the differential field $(C\langle
x,d\rangle^{(F^*,H)}, \delta )$.

Note that the field $C\langle x,d\rangle^{(F^*,H)}$ is not
invariant with respect to $d$.

The next result deals with generators of the differential field
 $(C\langle x,d\rangle^{(F^*,H)},\delta )$ over  $C$ and their relations. Let $p^d\langle x\rangle=\frac{
a^d\{x\}}{b^d\{x\}}$, where  $a^d\{x\},b^d\{x\}\in C\{x, d\}.$ Due
to (1) for any $g\in F^*$ one has $ga^d\{x\}b^{gd}\{x\}-
a^{gd}\{x\}b^d\{x\}= 0.$ Consider the left side of it as a
$d$-differential polynomial in $x$. In this case the above
equality means that all coefficients of this polynomial, which are
$d$-differential polynomials in $g$ over $C$, are zero for any
$g\in F^*$. Therefore due to assumption on $(F,d)$ for any
differential indeterminate $t$ one has $ta^d\{x\}b^{td}\{x\}-
a^{td}\{x\}b^d\{x\}= 0.$ Substitution $t= p^d\langle
x\rangle^{-1}$ implies that $p^{\delta}\langle x\rangle = 1 $

It should be noted that the field $C\langle x,d\rangle^H$ is
invariant with respect to the differential operator $d$ and
$(C\langle x,d\rangle^H,d)$ is a finitely generated
$d$-differential field over $C$ as a subfield of $C\langle x,d
\rangle.$ Let us assume that $(C\langle x,d \rangle^H, d)=
C\langle\varphi^d_1\langle x\rangle, \varphi^d_2\langle
x\rangle,...,\varphi^d_m\langle x\rangle,d \rangle.$ As far as
$p^d\langle x\rangle \in C\langle x,d \rangle^H$ there is a
differential rational function $\overline{p}^{d}\langle t_1, t_2
,...,t_m\rangle$ such that
 $$p^d\langle x\rangle =\overline{p}^{d}\langle \varphi^d_1\langle x\rangle, \varphi^d_2\langle x\rangle,...,\varphi^d_m\langle x\rangle\rangle.$$
Therefore due to $p^{\delta}\langle x\rangle= 1$ one has
$\overline{p}^{\delta}\langle \varphi^{\delta}_1\langle x\rangle,
\varphi^{\delta}_2\langle x\rangle,...,\varphi^{\delta}_m\langle
x\rangle\rangle =1. $

 {\bf Theorem 1.} \textit{ If $(C\langle x,d \rangle^H, d)= C\langle\varphi^d_1\langle x\rangle, \varphi^d_2\langle x\rangle,...,\varphi^d_m\langle x\rangle,d \rangle $ then the $\delta$- differential field $(C\langle x,d\rangle^{(F^*,H)},\delta )$ is generated over $C$ by the system of elements $\varphi^{\delta}_1\langle x\rangle, \varphi^{\delta}_2\langle x\rangle,...,\varphi^{\delta}_m\langle x\rangle$. Moreover any $\delta$-differential polynomial relation over $C$ of the system
 $\varphi^{\delta}_1\langle x\rangle, \varphi^{\delta}_2\langle x\rangle,...,\varphi^{\delta}_m\langle x\rangle $ is a
consequence of $d$-differential polynomial relations of the system
$\varphi^d_1\langle x\rangle, \varphi^d_2\langle
x\rangle,...,\varphi^d_m\langle x\rangle $ over $C$ and the
relation $\overline{p}^{\delta}\langle \varphi^{\delta}_1\langle
x\rangle, \varphi^{\delta}_2\langle
x\rangle,...,\varphi^{\delta}_m\langle x\rangle\rangle =1.$}

{\bf Proof.} It is evident that $C\langle x,d\rangle^{(F^*,H)}=
C\langle x,d\rangle^{F^*}\cap C\langle x,d\rangle^H.$

Let $ f^d\langle \varphi^d_1\langle x\rangle, \varphi^d_2\langle
x\rangle,...,\varphi^d_m\langle x\rangle\rangle$ be
$F^*$-invariant. Consider
$f^{t^{-1}d}\langle\varphi^{t^{-1}d}_1\langle x\rangle,
\varphi^{t^{-1}d}_2\langle x\rangle,...,\varphi^{t^{-1}d}_m\langle
x\rangle\rangle$ as a differential rational function in $x$ over
$C\langle t,d\rangle$ and let $\frac{P^d_t\{x\}}{Q^d_t\{x\}}$ be
its irreducible representation. $F^*$-invariantness of  $
f^d\langle \varphi^d_1\langle x\rangle, \varphi^d_2\langle
x\rangle,...,\varphi^d_n\langle x\rangle\rangle$ implies that
$\frac{P^d_t\{x\}}{Q^d_t\{x\}}= \frac{P^d_1\{x\}}{Q^d_1\{x\}}$. In
its turn it implies that
$$P^d_1\{x\}= P^d_t\{x\}\chi^d\langle t\rangle$$
, where $\chi^d\langle t\rangle \in C\langle t,d\rangle.$ Due to
Proposition 1 one has $\chi^d\langle p^d\langle x\rangle \rangle
\neq 0$ therefore  the above equality implies that
$$\frac{P^{d}_{p^d\langle x\rangle}\{x\}}{Q^d_{p^d\langle x\rangle}\{x\}}= \frac{P^d_1\{x\}}{Q^d_1\{x\}} \quad
\mbox{i.e.} \quad f^{\delta}\langle \varphi^{\delta}_1\langle
x\rangle,\varphi^{\delta}_2\langle
x\rangle,...,\varphi^{\delta}_m\langle  x\rangle\rangle=
f^d\langle \varphi^d_1\langle x\rangle, \varphi^d_2\langle
x\rangle,...,\varphi^d_m\langle x\rangle\rangle $$ So it is shown
that $C\langle x,d\rangle^{(F^*,H)}=C\langle
\varphi^{\delta}_1\langle x\rangle, \varphi^{\delta}_2\langle
x\rangle,...,\varphi^{\delta}_m\langle x\rangle ,\delta\rangle$.

Now let $\psi^{\delta}\{\varphi^{\delta}_1\langle x\rangle,
\varphi^{\delta}_2\langle x\rangle,...,\varphi^{\delta}_m\langle
x\rangle \} = 0$, where $\psi^d\{ t_1,t_2,...,t_m \} \in C\{
t_1,t_2,...,t_m,d\}$.

If $\psi^{d}\{\varphi^{d}_1\langle x\rangle, \varphi^{d}_2\langle
x\rangle,...,\varphi^{d}_m\langle x\rangle \} = 0$ then it means
that the above relation ($\psi^{\delta}\{ t_1,t_2,...,t_m \}$) of
the system $\varphi^{\delta}_1\langle x\rangle,
\varphi^{\delta}_2\langle x\rangle,...,\varphi^{\delta}_n\langle
x\rangle $ is a consequence of the relation ($\psi^{d}\{
t_1,t_2,...,t_m \}$) of the system $\varphi^d_1\langle x\rangle,
\varphi^d_2\langle x\rangle,...,\varphi^d_m\langle x\rangle $ i.e.
it is obtained by substitution $\delta$ for $d$ in $\psi^{d}\{
t_1,t_2,...,t_m \}$.

If $\psi^{d}\{\varphi^{d}_1\langle x\rangle, \varphi^{d}_2\langle
x\rangle,...,\varphi^{d}_m\langle x\rangle \} \neq 0$ then
consider $\psi^{t^{-1}d}\{\varphi^{t^{-1}d}_1\langle x\rangle,
\varphi^{t^{-1}d}_2\langle x\rangle,...,\varphi^{t^{-1}d}_m\langle
x\rangle \}$ as a $d$-differential rational function in one
variable $t$ over $C\langle x,d\rangle$. Let $\frac{a^d_x\{
t\}}{b^d_x\{ t\}}$ be its irreducible representation and the
leading coefficient (with respect to some linear order) of
$b^d_x\{t\}$ be one. We show that in this case all coefficients of
$a^d_x\{ t\}, b^d_x\{t\}$ belong to $C\langle x,d\rangle^H$.

Indeed, first of all $\psi^{t^{-1}d}\{\varphi^{t^{-1}d}_1\langle
x\rangle, \varphi^{t^{-1}d}_2\langle
x\rangle,...,\varphi^{t^{-1}d}_m\langle x\rangle \}$, as a
differential rational function in $x$, is $H$- invariant function,
as much as $\varphi^d_i\langle x\rangle \in C\langle
x,d\rangle^H$. This $H$-invariantness implies that
$$a^d_x\{t\}b^d_{hx+h_0}\{ t\}= b^d_x\{t\}a^d_{hx+h_0}\{t\}$$ for any $(h,h_0)\in H$. Therefore
$b^d_{hx+h_0}\{ t\}= \chi^d\langle x,(h,h_0)\rangle b^d_x\{t\}$.
But comparison of the leading terms of both sides implies that in
reality $\chi^d\langle x,(h,h_0)\rangle= 1$ which in its turn
implies that all coefficients of $b^d_x\{t\}$ (as well as
$a^d_x\{t\}$) are $H$- invariant.

Therefore all coefficients of $a^d_x\{t\}$, $b^d_x\{t\}$ can be
considered as $d$-differential rational functions in
$\varphi^d_1\langle x\rangle, \varphi^d_2\langle
x\rangle,...,\varphi^d_m\langle x\rangle$, for example
$b^d_x\{t\}= \overline{b}_{\varphi^d_1\langle x\rangle,
\varphi^d_2\langle x\rangle,...,\varphi^d_m\langle
x\rangle}\{t\}$. Now represent the numerator $a^d_x\{t\} $ as a
$d$-differential polynomial function in $t- \overline{p}^d\langle
\varphi^d_1\langle x\rangle, \varphi^d_2\langle
x\rangle,...,\varphi^d_m\langle x\rangle\rangle $, for example
$a^d_x\{t\}=\overline{a}^d_{\varphi^d_1\langle x\rangle,
\varphi^d_2\langle x\rangle,...,\varphi^d_m\langle x\rangle}\{t-
\overline{p}^d\langle \varphi^d_1\langle x\rangle,
\varphi^d_2\langle x\rangle,...,\varphi^d_m\langle x\rangle\rangle
\}$ . As such polynomial its constant term is zero because of
$\psi^{\delta}\{\varphi^{\delta}_1\langle x\rangle,
\varphi^{\delta}_2\langle x\rangle,...,\varphi^{\delta}_m\langle
x\rangle \} = 0$. So $$\psi^{t^{-1}d}\{\varphi^{t^{-1}d}_1\langle
x\rangle, \varphi^{t^{-1}d}_2\langle
x\rangle,...,\varphi^{t^{-1}d}_m\langle x\rangle \}
=\frac{\overline{ a}^d_{\varphi^d_1\langle x\rangle,
\varphi^d_2\langle x\rangle,...,\varphi^d_m\langle
x\rangle}\{t-\overline{p}^d\langle \varphi^d_1\langle x\rangle,
\varphi^d_2\langle x\rangle,...,\varphi^d_m\langle x\rangle
\rangle \}}{\overline{b}^d_{\varphi^d_1\langle x\rangle,
\varphi^d_2\langle x\rangle,...,\varphi^d_m\langle
x\rangle}\{t\}}.$$

Substitution $t=1$ implies that
$$\psi^{d}\{\varphi^{d}_1\langle x\rangle, \varphi^{d}_2\langle x\rangle,...,\varphi^{d}_m\langle x\rangle \} =\frac{\overline{
a}^d_{\varphi^d_1\langle x\rangle, \varphi^d_2\langle
x\rangle,...,\varphi^d_m\langle x\rangle}\{1-\overline{p}^d\langle
\varphi^d_1\langle x\rangle, \varphi^d_2\langle
x\rangle,...,\varphi^d_m\langle x\rangle \rangle
\}}{\overline{b}^d_{\varphi^d_1\langle x\rangle,
\varphi^d_2\langle x\rangle,...,\varphi^d_m\langle
x\rangle}\{1\}}.$$

Now consider the following $\delta$-differential rational function
over $C$:
$$\overline{\psi}^{\delta}\langle t_1,t_2,...,t_m \rangle=
\psi^{\delta}\{ t_1,t_2,...,t_m \}-\frac{\overline{
a}^{\delta}_{t_1, t_2,...,t_m}\{1-\overline{p}^{\delta}\langle
t_1,t_2,...,t_m \rangle \}}{b^{\delta}_{t_1,t_2,...,t_m}\{1\}}.$$
For this function one has
$$\overline{\psi}^{\delta}\langle \varphi^{\delta}_1\langle x\rangle, \varphi^{\delta}_2\langle x\rangle,...,\varphi^{\delta}_m\langle x\rangle \rangle= 0 \quad
\mbox{as well as} \quad \overline{\psi}^d\langle
\varphi^{d}_1\langle x\rangle, \varphi^{d}_2\langle
x\rangle,...,\varphi^{d}_m\langle x\rangle\rangle = 0.$$ But once
again the last equality(relation) means that it is consequence of
relations of the system $\varphi^{d}_1\langle x\rangle,
\varphi^{d}_2\langle x\rangle,...,\varphi^{d}_m\langle x\rangle $.
This is the end of proof of Theorem 1.

 {\bf Example 3.} Let $C=R$, $n=2$ and
$$H= O(2,R)=\{h\in GL(2,R): hh^t= E \}$$ In this case $R\langle x,
d \rangle^H= R\langle \varphi^d_1\langle x\rangle,
\varphi^d_2\langle x\rangle;d  \rangle$, where $\varphi^d_1\langle
x\rangle =(x,x), \varphi^d_2\langle x\rangle= (dx,dx)$ and $(*,*)$
stands for the dot product. Moreover $\varphi^d_1\langle x\rangle,
\varphi^d_2\langle x\rangle$ is $d-$algebraic independent over
$R$. For $p^d\langle x\rangle$ one can take $(x,dx)$, so $\delta=
(x,dx)^{-1}d$. Due to the Theorem 1
$$(R\langle x,d\rangle^{(F^*,H)},\delta )= R\langle
\varphi^{\delta}_1\langle x\rangle, \varphi^{\delta}_2\langle
x\rangle x\rangle ;\delta \rangle$$
In this case
$\overline{p}^d\langle \varphi^d_1\langle x\rangle
,\varphi^d_2\langle x\rangle \rangle=
\frac{1}{2}d\varphi^d_1\langle x\rangle$. So $p^{\delta}\langle
x\rangle= \frac{1}{2}\delta\varphi^{\delta}_1\langle x\rangle =
\frac{1}{2}\frac{1}{(x,dx)}d(x,x)=1$ and
$\varphi^{\delta}_2\langle x\rangle =(\delta x,\delta x)$ is
$\delta$-algebraic independent over $R$.

{\bf Example 4.} Let now $H= O(2,R)\propto R^2$. In this case
$R\langle x, d \rangle^H= R\langle \varphi^d_1\langle x\rangle,
\varphi^d_2\langle x\rangle;d  \rangle$, where $\varphi^d_1\langle
x\rangle =(dx,dx), \varphi^d_2\langle x\rangle= (d^2x,d^2x),$ and
$\varphi^d_1\langle x\rangle, \varphi^d_2\langle x\rangle$ is
$d-$algebraic independent over $R$. For $p^d\langle x\rangle$ one
can take
$$d\frac{\det[dx,d^2x]^2}{(dx,dx)^3}=d\frac{\varphi^d_1\langle
x\rangle\varphi^d_2\langle x\rangle-
\frac{1}{4}(d\varphi^d_1\langle x\rangle )^2}{(\varphi^d_1\langle
x\rangle)^3}= $$ $$\overline{p}^d\langle \varphi^d_1\langle
x\rangle ,\varphi^d_2\langle x\rangle \rangle $$ Due to the
Theorem 1
$$(R\langle x,d\rangle^{(F^*,H)},\delta )=
R\langle \varphi^{\delta}_1\langle x\rangle,
\varphi^{\delta}_2\langle x\rangle ;\delta \rangle$$ In this case
we have one basic relation $$\overline{p}^{\delta}\langle
\varphi^{\delta}_1\langle x\rangle ,\varphi^{\delta}_2\langle
x\rangle \rangle =1$$

The next important problem is the finding out the number of
elements of a maximal $\delta$-algebraic independent system of
elements of $C\langle x, d\rangle^{(F^*,H)}$ over $C$ i.e.
evaluation the differential transcendence degree
$\delta-\mbox{tr.deg} C\langle x, d\rangle^{(F^*,H)}/C $ of
$C\langle x, d\rangle^{(F^*,H)}$ over the its field of constants
$C$. Theorem 2 deals with this problem.

 {\bf Theorem 2.}  $\delta-\mbox{tr.deg} C\langle x, d\rangle^{(F^*,H)}/C = n-1$.

{\bf Proof.} Indeed due to Theorem 1, applied to H consisting of
only identity element, $C\langle x, d\rangle^{F^*}= C\langle x,
\delta \rangle$ with only one main relation among
$x_1,x_2,...,x_n$, namely $p^{\delta}\langle x\rangle = 1 $ i. e.
$a^{\delta}\{x\}-b^{\delta}\{x\}=0.$ But $a^d\{x\}-b^d\{x\}$ is
not zero polynomial as far as $p^{d}\langle x\rangle \neq 1 $.
Therefore  $\delta-\mbox{tr.deg} C\langle x,\delta\rangle/C =
n-1$.

Now to establish the equality $\delta-\mbox{tr.deg} C\langle x,
d\rangle^{(F^*,H)}/C = n-1$ it is enough to show that every $x_i$
is $\delta$-differential algebraic over $(C\langle
x\rangle^{(F^*,H)}, \delta )$. For it consider the following
$$W^{\delta}(\delta x_1,\delta x_2,...,\delta x_n)^{-1}W^{\delta}(\delta x_1,\delta x_2,...,\delta x_n,\delta y )=0 $$
as a differential equation in one variable $y$.

One can check easily that
$$ W^{\delta}(\delta x_1,\delta x_2,...,\delta x_n)^{-1}W^{\delta}(\delta x_1,\delta x_2,...,\delta x_n,\delta y )= \sum^{n+1}_{i=1}(-1)^{n+1-i}\frac{W^{\delta}_i}{W^{\delta}}\delta^iy $$
But $\frac{W^{\delta}_i}{W^{\delta}}\in C\langle x,
d\rangle^{(F^*,H)}$ and $y= x_i$ is a solution of the above
differential equation for any $i=\overline{1,n}$, which implies
that $\delta-\mbox{tr.deg} C\langle x, d\rangle^{(F^*,H)}/C =
n-1.$ This is the end of proof of Theorem 2.

{\bf 4. On differential field $(C\langle x,d\rangle^H,d)$}

Due to Theorem 1 one can find a system of differential generators
of $(C\langle x,d\rangle^H, d)$ to get a system of differential
generators of $(C\langle x, d\rangle^{(F^*,H)},\delta )$. In [2]
it is proved that $C\langle x,d\rangle^{GL(n,C)\propto C^n}$ as a
$d$-differential field over $C$ is generated by the following
$d$-algebraic independent system
$\frac{W^d_1}{W^d},\frac{W^d_2}{W^d},...,\frac{W^d_n}{W^d}$.

The next result reduces the problem of description of the
differential field $(C\langle x,d\rangle^H,d)$ to the description
of the ordinary field of algebraic invariants of $H$.

 {\bf Theorem 3.} \textit{The system $x,dx, d^2x,...,d^nx$ is algebraic independent over
 $C\langle x,d\rangle^{GL(n,C)\propto C^n}$ and
 $C\langle x,d\rangle=C\langle x,d\rangle^{GL(n,C)\propto C^n}(x,dx, d^2x,...,d^nx)$.}

{\bf Proof.} It is evident that it is enough to prove the
algebraic independence of the system $dx, d^2x,...,d^nx$  over
$C\langle x,d\rangle^{GL(n,C)\propto C^n}$. Let
$P[z_1,z_2,...,z_n]$, where $z_i= (z_{1i},z_{2i},...,z_{ni})$, be
a nonzero polynomial over $C\langle x,d\rangle^{GL(n,C)\propto
C^n}$ such that $P[dx,d^2x,...,d^nx]= 0$.  Assume, for example, at
least one of $z_{n1},z_{n2},...,z_{nn} $ occurs in $P$ and
$$P[z_1,z_2,...,z_n]= \sum_{\alpha}(z_{n1})^{\alpha_1}(z_{n2})^{\alpha_2}...(z_{nn})^{\alpha_n}P_{\alpha}[\overline{z_1},\overline{z_2},...,\overline{z_n}]$$, where $P_{\alpha}[\overline{z_1},\overline{z_2},...,\overline{z_n}]$ are polynomials over $C\langle x,d\rangle^{GL(n,C)\propto C^n}$ in  $\overline{z_i}=(z_{1i},z_{2i},...,z_{(n-1)i})$, $i=\overline{1,n}$. To prove Theorem 3 it is enough to show that $P[dx,d^2x,...,d^nx]= 0$ implies $P_{\alpha}[\overline{dx},\overline{d^2x},...,\overline{d^{n}x}]=0$ for all $\alpha$.

Consider $h \in GL(n,C)$ which's $i$-th column is of the form
$(0,...,0,1,0,...,0,c_i)$, where $i=\overline{1,n-1}$ and its
$n$-th column is $(0,...,0,c_n)$. For such $h$ one has
$\overline{d^i(hx)}= \overline{d^ix}$. So far as the coefficients
of $P[z_1,z_2,...,z_n]$ are $GL(n,C)\propto C^n$- invariant,
substitution $hx$ for $x$ into $P[dx,d^2x,...,d^nx]= 0$ implies
that
$$\sum_{\alpha}(\sum_{i=1}^nc_idx_i)^{\alpha_1}(\sum_{i=1}^nc_id^2x_i)^{\alpha_2}
...(\sum_{i=1}^nc_id^nx_i)^{\alpha_n}P_{\alpha}[\overline{dx},
\overline{d^2x},...,\overline{d^nx}]= 0$$ Therefore due to the
assumption on $(F,d)$ for variables $y_1,y_2,...,y_n$ one has
\begin{eqnarray}\begin{array}{c}
\sum_{\alpha}(\sum_{i=1}^ny_idx_i)^{\alpha_1}(\sum_{i=1}^ny_id^2x_i)^{\alpha_2}...(\sum_{i=1}^ny_id^nx_i)^{\alpha_n}P_{\alpha}[\overline{dx},\overline{d^2x},...,\overline{d^nx}]=
0 \end{array}\end{eqnarray}

Now consider the ring $C\langle x,d\rangle[y_1,y_2,...,y_n]$ with
respect to differential operators\\ $\partial_1=
\frac{\partial}{\partial y_1},\partial_2= \frac{\partial}{\partial
y_2}, ...,\partial_n= \frac{\partial}{\partial y_n}$. It is clear
that its constant ring is $C\langle x,d\rangle$ i.e. $$C\langle
x,d\rangle= \{a\in C\langle x,d\rangle[y_1,y_2,...,y_n]:
\partial_1a=\partial_2a=...\partial_na=0 \}.$$ Introduce new
differential operators $\overline{\partial_i}=
\sum_{j=1}^nf^d_{ij}\langle x \rangle\partial_j$, where
$i=\overline{1,n}$,$$(f^d_{ij}\langle x
\rangle)_{i,j=\overline{1,n}}=[dx,d^2x,...,d^nx]^{-1}$$ The
following are evident:

 a) The constant ring of $C\langle x,d\rangle[y_1,y_2,...,y_n]$ with respect to new differential operators $\overline{\partial_1},\overline{\partial_2},...,\overline{\partial_n}$ is the same $C\langle x,d\rangle$,

b) $\overline{\partial_j}(\sum_{i=1}^ny_id^kx_i)$ is equal to  $0$
whenever
 $j \neq k$ and it is equal to 1 if $j= k$.

Now if one assumes that $\alpha^0= (\alpha^0_1,...,\alpha^0_n)$ is
a such one for which $$\vert \alpha^0 \vert = \mbox{ max}\{
\vert\alpha \vert:
P_{\alpha}[\overline{dx},\overline{d^2x},...,\overline{d^nx}]\neq
0 \}$$ and applies
$(\overline{\partial_1})^{\alpha^0_1}(\overline{\partial_2})^{\alpha^0_2}...(\overline{\partial_n})^{\alpha^0_n}$
to equality (5) he comes to a contradiction
$P_{\alpha^0}[\overline{dx},\overline{d^2x},...,\overline{d^nx}]=
0$.

The equality $C\langle x, d\rangle= C\langle x,
d\rangle^{GL(n,C)\propto C^n}(x,dx,...,d^nx)$ is an easy
consequence of the fact that $y= x_i$, $i=\overline{1,n}$, are
solutions of the differential equation
$$W(dx_1,dx_2,...,dx_n)^{-1}W(dx_1,dx_2,...,dx_n,dy)= 0$$  all coefficients of which belong to $C\langle x, d\rangle^{GL(n,C)\propto C^n}$.
 This is the end of proof Theorem 3.

So due to Theorem 3 $$C\langle x, d\rangle^H= C\langle x,
d\rangle^{GL(n,C)\propto C^n}(x,dx,...,d^nx)^H$$ and the system
$\{x,dx,...,d^nx\}$ is algebraic independent over $C\langle x,
d\rangle^{GL(n,C)\propto C^n}$. But as we have already noticed the
$d$-differential field $C\langle x, d\rangle^{GL(n,C)\propto C^n}$
is generated over $C$ by the $d$-algebraic independent system
$\frac{W^d_1\{dx,d^{n+1}x\}}{W^d\{dx\}},\frac{W^d_2\{dx,d^{n+1}x\}}{W^d\{dx\}},...,\frac{W^d_n\{dx,d^{n+1}x\}}{W^d\{dx\}}$.
Therefore to find a system of generators of the field $C\langle x,
d\rangle^{GL(n,C)\propto C^n}(x,dx,...,d^nx)^H$ over $C\langle x,
d\rangle^{GL(n,C)\propto C^n}$  it is enough to find a system of
generators of the field $C(x,dx,...,d^nx)^H$ over $C$( see the
Appendix). Due to these the following method can be used to
construct a system of generators of the d-differential field
$C\langle x, d\rangle^H$ over $C$:

1) Find any system of ordinary algebraic generators of the field
$C(z_1,z_2,...,z_{n+1})^H$, where $z_i=(z_{1i},z_{2i},...,z_{ni})$
, $i=\overline{1,n+1}$, and the action of $H$ is defined as:
$$((h,h_0),(z_1,z_2,...,z_{n+1}))\rightarrow
(hz_1+h_0,hz_2,...,hz_{n+1})$$  For example, let it be
$\varphi_1(z_1,z_2,...,z_{n+1}),
\varphi_2(z_1,z_2,...,z_{n+1}),...,\varphi_k(z_1,z_2,...,z_{n+1})$.

2) The system
$$\frac{W^d_1\{dx,d^{n+1}x\}}{W^d\{dx\}},\frac{W^d_2\{dx,d^{n+1}x\}}{W^d\{dx\}},...,
\frac{W^d_n\{dx,d^{n+1}x\}}{W^d\{dx\}},\varphi_1(x,dx,...,d^nx),...,\varphi_k(x,dx,...,d^nx)$$
will be a system of generators of the differential field
$(C\langle x, d\rangle^H, d)$ over $C$.

{\bf 5. Discussion of two questions}

The following kind questions are typical in Differential Geometry.
The complete answers to these questions can not be given in our
case because of its generality.  The answers depend on
differential field $(F,d)$ and group $H$.

The first one is if $C\langle x, d\rangle^{(F^*,H)}$ distinguishes
nonequivalent curves ("of common position"). In our case it should
be noted that if pairs $(u,d_1)$, $(v, d_2)$ are
$(F^*,H)$-equivalent then $\delta(u,d_1)= \delta(v,d_2)$ i.e.
$(p^{d_1}\langle u\rangle)^{-1}d_1= (p^{d_2}\langle
v\rangle)^{-1}d_2$, provided that $p^{d_1}\langle u\rangle$,
$p^{d_2}\langle v\rangle$ are not zero. Therefore in our case the
equality $\delta_0= (p^{d_1}\langle u\rangle)^{-1}d_1=
(p^{d_2}\langle v\rangle)^{-1}d_2$ is a natural condition for the
pairs $(u,d_1)$, $(v, d_2)$ to be $(F^*,H)$-equivalent.

Consider $X= \{ v\in F^n: W^d\{dv\}p^d\langle v\rangle \neq 0 \}$,
which is not empty due to the assumption on differential field
$(F, d)$. Assume that $(u,d_1)$, $(v, d_2)$ are  such pairs that
$\delta_0= \delta(u,d_1)= \delta(v,d_2)$ and for every
$\varphi^{\delta}\langle x\rangle \in C\langle x,
d\rangle^{(F^*,H)}$ the equality $\varphi^{\delta_0}\langle
u\rangle = \varphi^{\delta_0}\langle v\rangle$ is valid. Then in
particular
$$a_i= \frac{W^{\delta_0}_i\{\delta_0u,\delta_0^{n+1}u\}}{W^{\delta_0}\{\delta_0u\}}=
\frac{W^{\delta_0}_i\{\delta_0v,\delta_0^{n+1}v\}}{W^{\delta_0}\{\delta_0v\}}
$$ at any $i=\overline{1,n}$. Which means that components both of $u$, $v$ are solutions of the same $(n+1)$-order linear differential equation
$$\delta_0^{n+1}y+ \sum_{i=1}^{n}(-1)^{n+1-i}a_i\delta_0^iy= 0 $$
one more solution of which is $y=1$. It means that $v=hu+ h_0$ for
some $(h,h_0)\in GL(n,C)\propto C^n$. Therefore in common case,
the problem is reduced to the question if
  $C(z_1,z_2,...,z_{n+1})^H $ distinguishes
$u,du,...,d^nu$ and $hu+h_0,hdu,...,hd^nu$.

 Another question is: What values can take a given system generators of  $(C\langle x, d\rangle^{(F^*,H)}, \delta)$?
Let us consider the case when the system of generators of
$(C\langle x, d\rangle^H,d)$ is given in the above considered
form: $$\frac{W^d_i\{dx,d^{n+1}x\}}{W^d\{dx\}},i=\overline{1n},
\qquad \varphi_j(x,dx,...,d^nx),j=\overline{1,k} $$ and let us
assume that $\{ \psi^d_i\{y_1,y_2,...,y_n,
t_1,t_2,...,t_k\}_{i=\overline{1,l}}$ is a fundamental system of
relations of the above system i.e. any other relation is an
element of the radical differential ideal generated by this system
in $C\{y_1,y_2,...,y_n, t_1,t_2,...,t_k; d\}$. Then due to Theorem
1 the system
$$\frac{W^{\delta}_i\{\delta x,\delta^{n+1}x\}}{W^{\delta}\{\delta x\}},i=\overline{1n}, \qquad \varphi_j(x,\delta x,...,\delta^nx),j=\overline{1,k}  $$
is a system of generators of the differential field $(C\langle x,
d\rangle^{(F^*,H)}, \delta)$ and the system consisting of
$\overline{p}^{\delta}\langle y_1,y_2,...,y_n,
t_1,t_2,...,t_k\rangle -1, \psi^{\delta}_i\{y_1,y_2,...,y_n,
t_1,t_2,...,t_k\}$, where $i=\overline{1,l}$, is a fundamental
system of relations of this system of generators.

Therefore if $(a_1,a_2,...,a_n, b_1,b_2,...,b_k)\in F^{n+k}$ is a
fixed element then the question if one can find such
$u=(u_1,u_2,...,u_n) \in F^n $ for which
$\frac{W^{\delta_0}_i\{\delta_0u,\delta_0^{n+1}u\}}{W^{\delta_0}\{\delta_0u\}}=
a_i$, where $i=\overline{1,n}$ and $\varphi_j(u,\delta_0
u,...,\delta_0^nu)=b_j$, where $j=\overline{1,k}$, $\delta_0=
p^d\langle u\rangle^{-1}d$, is equivalent to the existence of a
solution in $F^n$ of the following system:
\begin{eqnarray} \left\{ \begin{array}{l}
\psi^{\delta}_i\{a_1,a_2,...,a_n, b_1,b_2,...,b_k\}=0, \quad \mbox{where}  \quad i=\overline{1,l}\\
\overline{p}^{\delta}\langle a_1,a_2,...,a_n, b_1,b_2,...,b_k\rangle -1= 0\\
\delta^{n+1}x+ \sum_{i=1}^{n}(-1)^{n+1-i}a_i\delta^ix= 0\\
\varphi_j(x,\delta x,...,\delta^nx)= b_j, \quad  \mbox{where }  \quad j=\overline{1,k}\\
\end{array} \right. \nonumber \end{eqnarray}
and, of course, $\delta = p^d\langle x\rangle^{-1}d$.

{\bf 6. Appendix}

Let $C$ be any field of characteristic zero, $t$ and $z$, where
$z$ stands for column vector with coordinates $z_1,...,z_n$, are
indeterminates over $C$, $H$ be any subgroup of $ GL(n,C)\propto
C^n$. For brevity, let us denote an element of $H$ as $h$ and its
action on $z$ as $hz$. We consider  $C(t,z)$ with respect to
transformations $(h;(t,z))\rightarrow (t,hz)$.

{\bf Proposition 2.} \textit{The field $C(t,z)^H$ is generated by
$C(z)^H$ over $C(t)$.}

{\bf Proof.} Let $\frac{P[z]}{Q[z]}$ be any irreducible ratio from
$C(t,z)^H$. One can assume that coefficients of  $P[z], Q[z]$ are
from $C[t]$. For any $h\in H$ the equality
$\frac{P[z]}{Q[z]}=\frac{P[hz]}{Q[hz]} $ implies that
\begin{eqnarray}\begin{array}{c}P[hz]=
P[z]\chi_{h}[z]\end{array}\end{eqnarray} , where $\chi_{h}[z]\in
C(t)[z]$. In particular $P[h^{- 1}z]= P[z]\chi_{h^{-1}}[z]$ i.e.
$P[z]= P[hz]\chi_h[hz]$. Combining it with (6) implies that
$P[hz]= P[hz]\chi_h[z]\chi_{h^{-1}}[hz]$ i.e.
$\chi_h[z]\chi_{h^{-1}}[hz]= 1$. It shows that $\chi_h[z]$ does
not depend on $z$, i.e.  $\chi_h[z]\in C(t)$. Now consider the
leader $z^{\alpha}$ of $P[z]$, with respect to any linear order,
and apply $(\frac{\partial}{\partial z_1})^
{\alpha_1}...(\frac{\partial}{\partial z_n})^ {\alpha_n}$ to
equality (6). It shows that in reality $\chi_h[z]$ is of the form
$\chi_h[z]=\frac{\phi_h[t]}{a[t]}$,
$\phi_h[t]=\sum_{i=0}^{k}\phi_i(h)t^i $, where $k$ is a
nonnegative integer, $\phi_i(h)\in C(h)$, $a[t]\in C[t]$, the
leading coefficient of $a[t]$ is assumed to be 1. Due to (6) the
polynomial $\phi_h[t]$ has property
\begin{eqnarray}\begin{array}{c} \phi_{h_1}[t]\phi_{h_2}[t]= \phi_{h_1h_2}[t]a[t]\end{array}\end{eqnarray}
for any $h_1, h_2\in H$. In particular it implies that
 $$\phi_k(h_1)\phi_k(h_2)= \phi_k(h_1h_2)$$ and therefore
 if $\phi_k[h_1]\neq 0$ for some $h_1\in H$ then $\phi_k[h]\neq 0$ for any $h\in H$. Now due to (7) it is clear that
$\mbox{deg}\phi_h[t]= \mbox{deg} a[t]$ for any $h\in H$. But due
to (7) if some prime $p[t]\in C[t]$ divides $a[t]$ then it divides
$\phi_h[t]$ as well. It is not difficult to see now that $a[t]$
divides $\phi_h[t]$ and therefore $\phi_h[t]= \phi_k(h)a[t]$ i.e.
$\chi_h[z]=\frac{\phi_h[t]}{a[x]}=\phi_k(h).$ Represent $P[z]$
$(Q[z]$ )in the form $\sum P_i[z]t^i $  (respect.$\sum Q_i[z]t^i
$), where $P_i[z]\in C[z]$  (respect. $Q_i[z]\in C[z]$), rewrite
(6) in the following form $\sum P_i[hz]t^i =\sum
P_i[z]t^i\phi_k(h)$ (respect.$\sum Q_i[hz]t^i =\sum
Q_i[z]t^i\phi_k(h)$). So for every $i$ one has
$P_i[hz]=P_i[z]\phi_k(h)$ (respect.$Q_i[hz]=Q_i[z]\phi_k(h)$) and
therefore $\frac{P[z]}{Q[z]}\in C(t)C(z)^H$.This is the end of
proof.

\begin{center}{References}\end{center}

1. U.D. Bekbaev, On differential rational invariants of finite
subgroups of Affine group. Bulletin of Malaysian Mathematical
Society, 2005, Volume 28, N1, pp. 55-60.

2.  Dj.Khadjiev, Application of Invariant  Theory to Differential
Geometry of curves, FAN, Tashkent,1988,(Russian).

3. U.D. Bekbaev, On invariants of Curves, INSPEM, Technical report
No.1 December,2002

4. U.D. Bekbaev, On invariants of curves. Proceedings of
International Conference on  Research and Education in Mathematics
, 2-4 April 2003. At: Hotel Equatorial, Kuala Lumpur, Malaysia.
pp. 82-95. University Putra Malaysia Press.

5. E.R.Kolchin,  Differential Algebra and Algebraic Groups,
Academic Press, New York, 1973.

6.  U.D. Bekbaev, Once again on equivalence and invariants of
differential  equations \\ $y^{(n)}+ a_{n-1}y^{(n-1)}+ . . . +
a_0y= 0$, Uzbek  Mathematical Journal.3(1995)19-31(Russian).

\end{document}